# $q$-Line Search Scheme for Optimization Problem

**Suvra Kanti Chakraborty, Geetanjali Panda**
Departement of Mathematics, Indian Insititute of Technology Kharagpur.

*Corresponding Author:*

First Author,
Departement of Mathematics,
Indian Institute of Technology Kharagpur,
Kharagpur, West Bengal, India, 721 302.
Email: suvrakanti@maths.iitkgp.ernet.in

**ABSTRACT**

In this paper new descent line search iterative schemes for unconstrained as well as constrained optimization problems are developed using $q$-derivative. At every iteration of the scheme, a positive definite matrix is provided which is neither exact Hessian of the objective function as in Newton scheme nor the positive definite matrix as generated in quasi-Newton scheme. Second order differentiability property is not required in this process. Component of this matrix are constructed using $q$-derivative of the function. It is proved that the schemes preserve the property of Newton-like schemes in a local neighborhood of a minimum point which leads to the super linear rate of convergence. Numerical illustration of the scheme is also provided.

*Keyword:*

$q$-derivative,

Newton like method,

unconstrained optimization,

SQP method.

## 1. INTRODUCTION

Since last few decades, $q$-calculus (quantum calculus) has been one of the interesting topics for both Mathematicians and Physicists. The $q$-analogue of ordinary derivative was first introduced by F. H. Jackson. Literature on the application of $q$-calculus is available in several areas (see [2], [9], [13],[14], [17], [19]). Some recent developments using $q$-derivatives can be found in transform calculus [1], variational calculus [3], sampling theory [15], $q$-version of Bochner Theorem [12], and so on. Soterroni et al. [20] first studied the use of $q$-derivative in the area of unconstrained optimization, which is the $q$-variant of steepest descent method. Gouvea et al. recently searched for global optimum [11] using $q$-steepest descent and $q$-CG method, where they have proposed a descent scheme using $q$-calculus with stochastic approach which does not focus order of convergence of the scheme. This paper has focused on a modified Newton-like deterministic scheme based on $q$-derivative which has super linear convergence property.

The general iterative line search process to minimize a function $f: R^n \to R$ is of the form as $x^{(k+1)} = x^{(k)} + \alpha_k p_k$, where $p_k$ is a descent direction at $x^{(k)}$ and $\alpha_k$ is the step length along $p_k$. A natural choice of $p_k$ is $p_k = -B_k^{-1} \nabla f(x^{(k)})$, where $B_k$ is a positive definite matrix. For Steepest Descent method and Newton method $B_k$ is taken to be an Identity matrix and exact Hessian of $f$ respectively. The quasi Newton BFGS scheme generates generates $B_k$ as

$$B_{k+1} = B_k - \frac{B_k s_k s_k^T B_k}{s_k^T B_k s_k} + \frac{y_k y_k^T}{y_k^T s_k},$$

where $s_k = x^{(k)} - x^{(k-1)}$ and $y_k = \nabla f(x^{(k)}) - \nabla f(x^{(k-1)})$. A good survey on line search schemes in this direction can be found in the book [16]. In this paper we propose a line search scheme which provides a positive definite matrix $B_k$ at $k^{th}$ iteration. $B_k$ is not generated from $B_{k-1}$ like quasi Newton schemes. For most of the successful Newton like schemes second order differentiability of the function is essential. Moreover, $B_k$ does not need second order partial derivatives of the function always and for large $k$, $B_k$ behaves like Hessian matrix of $f(x)$.

In Section 2, some notations and definitions on $q$-calculus and other prerequisites are provided, which are used in sequel throughout the paper. The scheme is proposed in Section 3 for unconstrained optimization problem. The convergence and advantages of the scheme is discussed in this section. The theoretical concept of the proposed scheme is extended to SQP scheme in Section 4. Section 5 concludes final remarks.

## 2. Prerequisites

For a function $f: \mathbb{R} \to \mathbb{R}$, the $q$-derivative ($q \neq 1$) of $f$ (denoted by $D_{q,x}f$), is defined as

$$D_{q,x}f(x) = \begin{cases} \dfrac{f(x) - f(qx)}{(1-q)x}, & x \neq 0 \\ f'(x), & x = 0 \end{cases}$$

Suppose $f: \mathbb{R}^n \to \mathbb{R}$, whose partial derivatives exist. For $x \in \mathbb{R}^n$, consider an operator $\epsilon_{q,i}$ on $f$ as

$$(\epsilon_{q,i}f)(x) = f(x_1, x_2, \ldots, x_{i-1}, qx_i, x_{i+1}, \ldots, x_n).$$

The $q$-partial derivative ($q \neq 1$) of $f$ at $x$ with respect to $x_i$, denoted by $D_{q,x_i}f$, is

$$D_{q,x_i}f(x) = \begin{cases} \dfrac{f(x) - (\epsilon_{q,i}f)(x)}{(1-q)x_i}, & x_i \neq 0 \\ \dfrac{\partial f}{\partial x_i}, & x_i = 0 \end{cases}$$

**Example:** $f(x, y) = y^2 + 4x^3$. Then $D_{q,x}f = 4x^2(1 + q + q^2)$, $D_{q,y}f = y(1 + q)$.

**Theorem 2.1.** *(Zoutendjik Theorem [16])*

*Consider $k^{th}$ iteration of an optimization algorithm in the form $x^{(k+1)} = x^{(k)} + \alpha_k p_k$, where $p_k$ is a descent direction and $\alpha_k$ satisfies Wolfe condition. Suppose $f$ is bounded below in $\mathbb{R}^n$ and that $f$ is continuously differentiable in an open set containing the level set $\mathcal{L} = \{x : f(x) \leq f(x_0)\}$, where $x_0$ is the starting point of the iteration.*

*Assume also that $\nabla f$ is Lipschitz continuous on $\mathcal{L}$. That is, there exists a constant $L_c > 0$ such that*

$$\|\nabla f(x) - \nabla f(\bar{x})\| < L_c \|x - \bar{x}\| \; \forall x, \bar{x} \in \mathcal{L},$$

*then*

$$\sum_{k \geq 0} \cos^2 \theta_k \left\|\nabla f(x^{(k)})\right\|^2 < \infty,$$

where $\theta_k$ is the angle between $p_k$ and $\nabla f(x^{(k)})$.

**Symmetric Indefinite Factorization([16]):**

A real symmetric matrix $A$ can be expressed as $PAP^T = LBL^T$, where $L$ is a lower triangular matrix, $P$ is a permutation matrix and $B$ is a block diagonal matrix which allows at most $2 \times 2$ blocks. This requires a pivot block initially. There are several pivoting strategies available in the literature (see [5], [7], [8], [10]) to take care the sparsity of the matrix. The symmetric indefinite factorization allows to determine the inertia of $B$ and inertia of $B$ remains equal to inertia of $A$. An indefinite factorization can be modified to ensure that the modified factors are the factors of a positive definite matrix. Consider the spectral decomposition of $B$ as $B = Q\Lambda Q^T$, where $Q$ is the matrix whose columns consist of eigen vectors and $\Lambda$ is the diagonal matrix whose diagonal elements are respective eigen values $\lambda_i$ of $B$. Denote $F = Q\, diag(\tau_i)\, Q^T$, where

$$\tau_i = \begin{cases} 0, & if\ \lambda_i \geq \delta \\ \delta - \lambda_i, & if\ \lambda_i < \delta \end{cases}$$

From the construction of $F$ it is true that $LBL^T$ is sufficiently positive definite. A matrix $E$ has to be added to $A$ to make it positive definite. $P(A+E)P^T = L(B+F)L^T$ provides $E = P^T LF L^T P$. So, $\lambda_{min}(A+E) \approx \delta$. In this process $\bar{A} \triangleq A + E$ is the positive definite matrix. This idea is briefed in [16] and for this purpose MATLAB in-built command $ldl(\ )$ is used in this paper since it is less expensive.

### 3. Proposing a New Scheme for Unconstrained Optimization Problem

Consider a general unconstrained optimization problem $(P): min_{x \in \mathbb{R}^n} f(x)$, where $q$-partial derivatives of $\nabla f$ exist. Consider a positive real sequence $\{q_k\}$ so that $q_k \to 1$ as $k \to \infty$. Some choices of $\{q_k\}$ may be taken as $q_{k+1} = 1 - \frac{q_k^\gamma}{k}$, with $0 < q_0 < 1$, $\gamma$ as a positive integer. At $x^{(k)}$, consider a matrix $A_{q_k} = (\frac{1}{2}(a_{ij} + a_{ji}))_{n \times n}$, where

$$a_{ij} = D_{q_k, x_i} \frac{\partial f(x^{(k)})}{\partial x_j}, \qquad 1 \leq i,j \leq n. \tag{1}$$

$A_{q_k}$ may not be positive definite, which can be modified to a positive definite matrix $B_{q_k}$ through symmetric indefinite factorization as described in previous section. With this modification the new iterative scheme becomes

$$x^{(k+1)} = x^{(k)} - \alpha_k B_{q_k}^{-1} \nabla f(x^{(k)}), \tag{2}$$

where $\alpha_k$ is the step length at $x^{(k)}$ along $-B_{q_k}^{-1} \nabla f(x^{(k)})$. Following algorithm describes the modified Newton like line search scheme using $q$-calculus.

Algorithm 1: $q$-line search for unconstrained optimization

Data: starting point $x^{(0)}$, $\epsilon$, $q_0$, $\gamma$;

for $k = 0, 1, 2 ...$
1. Compute $B_{q_k}$;
2. $x^{(k+1)} = x^{(k)} - \alpha_k B_{q_k}^{-1} \nabla f(x^{(k)})$, $\alpha_k$ satisfies Wolfe conditions;
3. if $\|\nabla f(x^{(k+1)})\| < \epsilon$
   Stop;

else

$k = k + 1, q_{k+1} = 1 - \frac{q_k^\gamma}{k}$;

end;

end;

## 3.1 Convergence of the proposed scheme

Following lemma guarantees the convergence of this new scheme.

**Lemma 3.1.** *Let $\kappa(B_{q_k})$ denotes the condition number of $B_{q_k}$. If there exists some $C > 0$ such that $\kappa(B_{q_k}) < C$ for every $k$, then under all the standard assumption of Zoutendjik Theorem 2.1, $\| \nabla f(x_k) \| \to 0$ as $k \to \infty$.*

**Proof**: Let the eigenvalues of $B_{q_k}$ be $0 < \lambda_1^{(k)} \leq \lambda_2^{(k)} \leq ... \leq \lambda_n^{(k)}$. Since $\lambda_1^{(k)}$ is the smallest eigenvalue of $B_{q_k}$, for any $u \in \mathbb{R}^n$,

$$u^T B_{q_k} u \geq \lambda_1^{(k)} \|u\|^2.$$

Let $\theta_k$ be the angle between $p_k$ and $\nabla f(x^{(k)})$, where $p_k = -B_{q_k}^{-1} \nabla f(x^{(k)})$. Hence,

$$\cos \theta_k = -\frac{\nabla f^{(k)T} p_k}{\|\nabla f^{(k)T}\| \|p_k\|} = \frac{p_k^T B_{q_k} p_k}{\|\nabla f^{(k)T}\| \|p_k\|} \geq \lambda_1^{(k)} \frac{\|p_k\|}{\|\nabla f(x^{(k)})\|}. \quad (3)$$

$\|\nabla f(x^{(k)})\| = \|B_{q_k} p_k\| \leq \|B_{q_k}\| \|p_k\| = \lambda_n^{(k)} \|p_k\|$. Using this in Eq. (3), we have

$$\cos \theta_k = -\frac{\nabla f^{(k)T} p_k}{\|\nabla f^{(k)T}\| \|p_k\|} \geq \frac{\lambda_1^{(k)}}{\lambda_n^{(k)}} = \frac{1}{\|B_{q_k}\| \|B_{q_k}^{-1}\|} \geq \frac{1}{C}.$$

Hence, under the assumption of Zoutendjik condition $\lim_{k \to \infty} \|\nabla f(x^{(k)})\| = 0$. □

In the closer neighborhood of the solution (if it is an isolated minimum point), the objective function shows strict convexity nature. As $k \to \infty$, $q_k \to 1$. Hence after some large $k$ it is obvious that the matrix $B_{q_k}$ becomes almost same as the exact Hessian of the objective function. So, the construction of different $\{q_k\}$ will certainly not affect the order of convergence of the scheme. Following theorem discusses order convergence.

**Theorem 3.1.** *Suppose that $f: \mathbb{R}^n \to \mathbb{R}$ is twice differentiable in the neighborhood of the local minimum point $x^*$ and the iterative scheme (2) converges to $x^*$. Then $\{x_k\}$ converges superlinearly iff $\lim_{k \to \infty} \frac{\|(B_{q_k} - \nabla^2 f(x^*))p_k\|}{\|p_k\|} = 0$.*

**Proof:** For large $k$, $\| \nabla^2 f(x^*) - B_{q_k} \| \leq \| \nabla^2 f(x^*) - \nabla^2 f(x^{(k)}) \| + \| B_{q_k} - \nabla^2 f(x^{(k)}) \| \to 0$. Hence by Theorem 3.7 of [16], the super linear rate of convergence of the proposed scheme is justified.

Some important features of the proposed scheme are:

- The scheme does not require second order differentiability like Newton method. It takes second order partial derivatives at most in a measure-zero set.
- The scheme provides positive definite matrix $B_{q_k}$ at each iteration, which is not generated from the matrix $B_{q_{k-1}}$ at previous iteration like Quasi-Newton method. $B_{q_k}$ uses value of $q_k$.
- The scheme has super linear convergence property.

### 3.2 Numerical Illustration

Following numerical illustration discusses the advantage of proposed $q$-line search method over classical Newton method. Consider a class of function generated by different values of $c \in \mathbb{R} \setminus \{0\}$.

$$f_c(x, y) = \begin{cases} 0.05(y - x^2)^2 + (1 - x)^2 + c, & x \geq c, \\ \frac{x}{c}(1 - x)^2 + 0.05(y - x^2)^2 - \frac{(1 - c)^2}{c}(x - c) + c, & x < c. \end{cases} \quad (4)$$

$f_c$, for $c \in \mathbb{R} \setminus \{0\}$ is not second order differentiable, since at $x = c$, $\frac{\partial^2 f_c}{\partial x^2}$ does not exist. It is also observed that each $f_c$ has minimum at $(1, 1)$ with minimum value $c$. So, classical Newton method can not be applied. But in this case $q$-line search method can be applied efficiently. For this purpose we have chosen some values of $c$ and constructed its corresponding $f_c$. To apply $q$-line search on $f_c$, we choose the initial guesses on $x = c$ line, which are of the form $(c, y)$ with $y = 0.1(0.2)1.9$.

For Algorithm 1, the terminating gradient norm is set as $10^{-5}$. 10 initial guesses are chosen on each of 10 lines given by $x = c$. The algorithm is executed in MATLAB R2015a software. We note the number of iterations and time elapsed in Table 1. So, we consider total 100 initial guesses which are inside the rectangle formed by the lines $x = 0.1$, $x = 1.9$, $y = 0.1$ and $y = 1.9$, which encloses $(1, 1)$. Fix $q_0$ as 0.9, as this choice of $q$ was previously taken in [18]. The step length $\alpha$ is chosen such that it satisfies Wolfe conditions with parameters $10^{-4}$ and 0.9 respectively. We compute the average iteration and average time elapsed and put a note of it for corresponding $f_c$. Besides the advantage of the proposed scheme over classical Newton method, we want to study a numerical comparison with existing first order method. Since BFGS scheme is treated to be one of the most efficient first order methods, we execute BFGS algorithm for each $f_c$ under the same environment with same initial guesses and other parameters.

In the following table for each chosen $c$, we note down the average iteration and average time required for both of the BFGS scheme and proposed method to find the minimum of corresponding $f_c$ taking the initial guesses as $(c, 0.1), (c, 0.3), \ldots, (c, 1.9)$. For this purpose we have chosen three variations $Q_\gamma$, $\gamma = 1, 2, 3$ of the proposed scheme formed by three sequences of $\{q_k\} \to 1$ as $q_{k+1} = 1 - \frac{q_k^\gamma}{k}$ respectively. We provided a graphical illustration in Figure 1, to see that the proposed method performs better with higher values of $\gamma$. One may also conclude from this numerical illustration that the proposed method solves the problem in lesser number of iterations and reasonable competitive time with BFGS method.

Table 1: Average Iteration and Time for $f_c$.

| c | BFGS | $Q_1$ | $Q_2$ | $Q_3$ | BFGS | $Q_1$ | $Q_2$ | $Q_3$ |
| --- | --- | --- | --- | --- | --- | --- | --- | --- |
| 0.1 | 7.1 | 5 | 5 | 4.9 | 3.308712431 | 4.1675119 | 4.1433847 | 4.0781232 |
| 0.3 | 7.3 | 5 | 4.9 | 4.7 | 3.377407186 | 4.1534731 | 4.1027449 | 3.9259322 |
| 0.5 | 9.8 | 5 | 4.8 | 4.5 | 4.264493988 | 4.2260957 | 4.0190671 | 3.7600182 |
| 0.7 | 8.1 | 4.6 | 4 | 4 | 3.742546297 | 3.9419581 | 3.4576616 | 3.4598311 |
| 0.9 | 7.5 | 4.1 | 3.7 | 3.3 | 3.505550687 | 3.5283773 | 3.2377628 | 2.9018503 |
| 1.1 | 8.0 | 4.1 | 3.8 | 3.7 | 3.689159394 | 3.705838 | 3.4486235 | 3.3703194 |
| 1.3 | 9.1 | 4.3 | 4.1 | 4 | 4.193454478 | 4.1288782 | 3.7385843 | 3.6639808 |
| 1.5 | 9.1 | 5.3 | 4.7 | 4.7 | 3.998260923 | 4.6753996 | 4.1988353 | 4.2017754 |
| 1.7 | 9.2 | 5.8 | 5.8 | 5.5 | 4.177781814 | 5.4869421 | 5.4943635 | 5.2393932 |
| 1.9 | 9.8 | 5.8 | 5.6 | 5.5 | 4.323570395 | 5.4316634 | 5.2661416 | 5.1906454 |

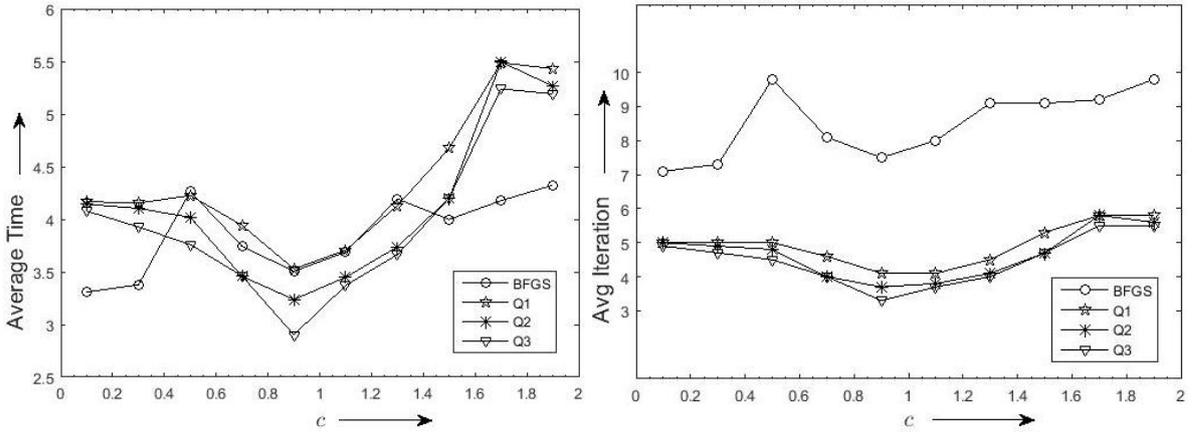

(a) Average Iteration for $f_c$

(b) Average Time for $f_c$

Figure 1: Comparison of Avg. Iteration and Time for $f_c$.

### 3.2.1 Test Problems

For numerical experiments, we measure the performance of BFGS algorithm and the proposed algorithm for a set of 15 problems in low dimensions taken from [21], some convex and some non convex. We run the codes for each problem in MATLAB R2015a platform and focus on iteration and time elapsed under the same environment. The name, dimension and minimum point are specified in Table 2. We follow the idea of Dolan and More [6] for preparing the performance profile. The final gradient norm for termination is taken as $10^{-5}$. We also stop if the time for a specific function exceeds 100 seconds. The step length $\alpha$ is chosen such that it satisfies Wolfe conditions with parameters $10^{-4}$ and 0.9. For practical purpose such $\alpha$ is found by backtracking with initial value 1 and backtracking factor $\tau = 0.5$. For the proposed scheme $q_0$ is chosen to be 0.9 and we take three variations of our schemes $Q_1, Q_2, Q_3$ as described in previous subsection. Since the schemes may lead to different minima of non convex function, we consider a run as a success if it takes to the global minimum point. The initial guesses are chosen randomly from a hypercube of length 1 keeping the desired minimum at the center of the hypercube. For preparing the performance profile, 10 successful runs are fixed for each scheme on each problem. Average iteration and average time in each case are noted.

The performance ratio defined by Dolan and More is $\rho_{(p,s)} = \frac{r_{(p,s)}}{\min\{r_{(p,s)}: 1 \leq r \leq n_s\}}$, where $r_{(p,s)}$ refers to the iteration (or, time) for solver $s$ spent on problem $p$ and $n_s$ refers to the number of problems in the model test set. In order to obtain an overall assessment of a solver on the given model test set, the cumulative distribution function $P_s(\tau)$ is $P_s(\tau) = \frac{1}{n_p} size\{p \in: \rho_{(p,s)} \leq \tau\}$, where $P_s(\tau)$ is the probability that a performance ratio $\rho_{(p,s)}$ is within a factor of $\tau$ of the best possible ratio.

Table 2: Name and Dimension of the function considered in Test set

| SlNo. | Function | Dim | Minimum point |
|---|---|---|---|
| 1 | $Bohachevsky$ | 2 | (0, 0) |
| 2 | $Branin$ | 2 | $(\pi, 2.275)$ |
| 3 | $Cross - in - Tray$ | 2 | (1.3491, 1.3491) |
| 4 | $Dixon - Price$ | 2 | (1, 0.7071) |
| 5 | $Easom$ | 2 | $(\pi, \pi)$ |
| 6 | $Griewank$ | 4 | (0, 0, 0, 0) |
| 7 | $Hartmann\ 3D$ | 3 | (0.1146, 0.5556, 0.8525) |
| 8 | $Levy$ | 4 | (1,1,1,1) |
| 9 | $Mccorm$ | 2 | $(-.54719, -1.54719)$ |
| 10 | $RotatedHyperEllipsoid$ | 4 | (0,0,0,0) |
| 11 | $Schwefel$ | 2 | (420.9687, 420.9687) |
| 12 | $Sphere$ | 8 | (0,0,0,0) |
| 13 | $Stybtang$ | 4 | $(-2.0953, -2.0953, -2.0953, -2.0953)$ |
| 14 | $Sum - Square$ | 10 | (0,0,,0) |
| 15 | $Zakharov$ | 2 | (0,0) |

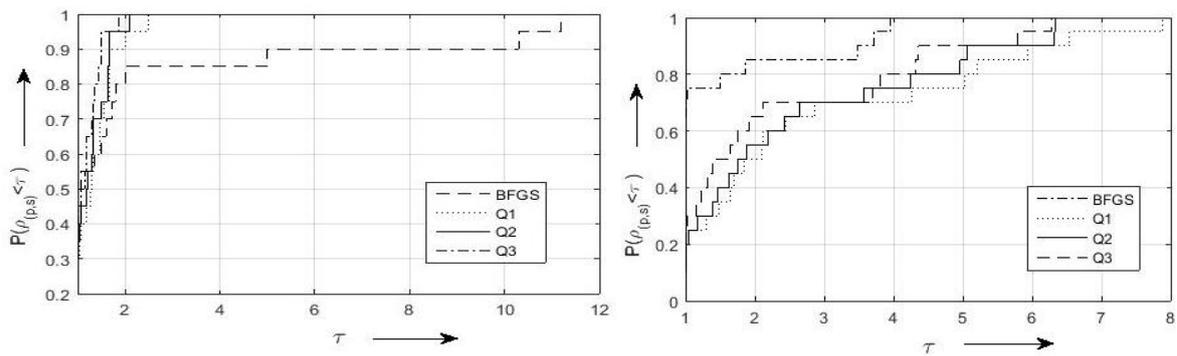

(a) Performance Profile of Iteration on Test Set $T$    (b) Performance Profile of Time on Test Set $T$

Figure 2: Performance Profile on Test Set $T$.

**4 Extension to Sequential Quadratic Programming(SQP)**

Consider the following constrained optimization problem:

$$(CP): \min_x f(x) \text{ subject to: } h(x) = 0, g(x) \leq 0,$$

where $f: \mathbb{R}^n \to \mathbb{R}$, $h: \mathbb{R}^n \to \mathbb{R}^m$ and $h: \mathbb{R}^n \to \mathbb{R}^p$. SQP is considered to be one of the most efficient method to solve (CP) for hassle free selection of starting point and its rapid convergence. The main idea of SQP ([16],[4]) is to model (CP) at $x^{(k)}$, by a quadratic subproblem QP, whose solution provides better approximation $x^{(k+1)}$. Finally the scheme reaches at the solution point $x^*$.

$$(QP): \min_{d_x} \nabla f(x^{(k)})^T d_x + \frac{1}{2} d_x^T B_k d_x$$

$$\text{subject to: } \nabla h(x^{(k)})^T d_x + h(x^{(k)}) = 0,$$

$$\nabla g(x^{(k)})^T d_x + g(x^{(k)}) \leq 0,$$

where $B_k$ is an approximate matrix of $\nabla_{xx}^2 \mathcal{L}(x^{(k)}, u^{(k)}, v^{(k)})$, and $\mathcal{L}$, $(u^{(k)}, v^{(k)})$ are corresponding Lagrangian function and Lagrange multipliers respectively. For Newton-SQP method, $B_k$ is exactly equal to $\nabla_{xx}^2 \mathcal{L}(x^{(k)}, u^{(k)}, v^{(k)})$ whereas for remote starting point $B_k$ can be chosen as positive definite approximation of $\nabla_{xx}^2 \mathcal{L}(x^{(k)}, u^{(k)}, v^{(k)})$. The positive definiteness of $B_k$ ensures the existence of unique solution of (QP). The updates of $(x, u, v)$ for some new step length $\alpha$ is

$$x^{(k+1)} = x^{(k)} + \alpha d_x, u^{(k+1)} = u^{(k)} + \alpha d_u, v^{(k+1)} = v^{(k)} + \alpha d_v.$$

Once the new iterates are constructed, the positive definite matrix $B_{k+1}$ is computed using the information of $B_k$. Local convergence proof is based on classical Newton method, which assumes $\alpha = 1$. But for remote starting point, $\alpha$ can be found through a merit function so that $\phi(x^{(k)} + \alpha d_x) < \phi(x^{(k)})$.

Instead of generating $B_k$ from $B_{k-1}$, we construct a positive definite matrix $B_{q_k}^L$ using the logic of Subsection 2 as follows. Consider a matrix $A_{q_k}^L = \left(\frac{a_{ij}^L + a_{ji}^L}{2}\right)_{n \times n}$, where $a_{ij}^L = D_{q_k, x_i} \frac{\partial L}{\partial x_j}(x^{(k)}, u^{(k)}, v^{(k)})$, $1 \leq i, j \leq n$. $B_{q_k}^L$ be the symmetric positive definite modification of $A_{q_k}^L$. The modified algorithm can be proposed as follows.

### 4.1 Convergence:

Let the approximation matrix $B_{q_k}^L$ satisfies the following assumptions.

**B1:** There exists a $\beta_1 > 0$ such that for each $k$, $d^T B_{q_k}^L d \geq \beta_1 \|d\|^2$ for all $d$ satisfying $\nabla h(x^{(k)})^T d = 0$.

**B2:** There exists a $\beta_2 > 0$ such that for each $k$, $\|B_{q_k}^L\| \leq \beta_2$.

**B3:** There exists a $\beta_3 > 0$ such that for each $k$, $B_{q_k}^{L^{-1}}$ exists and $\|B_{q_k}^{L^{-1}}\| \leq \beta_3$.

The convergence properties of SQP method depend on $B_{q_k}^L$ matrix. The rate of convergence of $x$-iterates can be achieved by approximating $B_{q_k}^L$ to $\nabla_{xx}^2 L(x^*, u^*, v^*)$ in an appropriate manner. The descent-ness and proper choice of $\alpha$ are necessary for global convergence of SQP schemes. For the discussion of rate of convergence we assume that the active inequality constraints for (CP) at $x^*$ are known. This assumption is quite justified because the problem (QP) at $x^{(k)}$ have the same active constraints as (CP) at $x^*$ when $x^{(k)}$ is near $x^*$. Hence the inequality constraints that are inactive at $x^*$ can be ignored and those that are active can be changed to equality constraints without changing the solution of (QP). So, under this assumption, equality constrained problem is considered for the local analysis. In this situation all dual variables are in the form of $u^{(k)}$ and $v^{(k)}$ is not reflected in this process. Now we recall the following theorem from Boggs and Tolle [4] for the optimization problem in the form of

$$\min_x f(x) \text{ subject to: } h(x) = 0.$$

**Theorem 4.1.** *[4] Let assumptions (B1-B3) hold and let the sequence $\{(x^{(k)}, u^{(k)})\}$ be generated by the SQP algorithm. Assume $x^{(k)} \to x^*$ and $\mathcal{P}$ be the projection operator on null space of $\nabla h^T$. Then the sequence $x^{(k)}$ converges to $x^*$ superlinearly iff the matrix approximations satisfy*

$$\lim_{k \to \infty} \frac{\|\mathcal{P}^K(B_k - \nabla_{xx}^2 \mathcal{L}((x^*, u^*)))(x^{(k+1)} - x^{(k)})\|}{\|(x^{(k+1)} - x^{(k)})\|} = 0. \quad (5)$$

*If this equation holds then the sequence $\{u^{(k)}\}$ converges R-superlinearly to $u^*$ and the sequence $\{(x^{(k)}, u^{(k)})\}$ converges superlinearly.*

---

**Algorithm 2:** $q$-line search SQP scheme for constrained optimization

---

Data: Input $(x^{(0)}, u^{(0)}, v^{(0)})$, $\epsilon$, $q_0$, $\gamma$;

for $k = 0, 1, 2 \ldots$

1. Compute $B_{q_k}^L$;
2. Solve (QP) with $(x^{(k)}, u^{(k)}, v^{(k)})$, and $B_{q_k}^L$ to obtain $(d_x, d_u, d_v)$;
3. Choose the step length $\alpha$ using a merit function;
4. Set

$x^{(k+1)} = x^{(k)} + \alpha d_x, u^{(k+1)} = u^{(k)} + \alpha d_u, v^{(k+1)} = v^{(k)} + \alpha d_v$;

5. Stop if $\|\nabla L\| < \epsilon$
6. Set $k = k + 1$, $q_{k+1} = 1 - \frac{q_k^\gamma}{k}$;

end;

---

In Algorithm 2, for large $k$, $B_k = B_{q_k}^L$ approximates $\nabla_{xx}^2 \mathcal{L}^k$ nicely. That is for arbitrary $\epsilon_1 > 0$, $\|B_{q_k}^L - \nabla_{xx}^2 \mathcal{L}^k\| < \epsilon_1$. Again, for arbitrary $\epsilon_2 > 0$, $\|\nabla_{xx}^2 \mathcal{L}^* - \nabla_{xx}^2 \mathcal{L}^k\| < \epsilon_2$. So, $\|(B_{q_k}^L - \nabla_{xx}^2 \mathcal{L}(x^*, u^*))\| \to 0$ as $k \to \infty$. So, by Theorem 4.1, the convergence rate are justified for the algorithm of this new $q$-line search SQP scheme.

Numerical experiments for (CP) with this new scheme can be done following the process of the numerical illustration of (P) in Subsection 3.2.

## 5 Conclusions

In this paper a new line search scheme is developed for unconstrained and constrained optimization problem. This scheme does not require second order differentiability assumption and the positive definite matrix at each iteration is free from the information of the positive definite matrix of previous iteration like general quasi Newton schemes. From the numerical illustrations it is observed that this new scheme can be a competitor of BFGS scheme under some favorable situations.